\documentclass{gtart}
\usepackage{amssymb,amsmath}
\usepackage{graphicx}

\input gtmonout
\volumenumber{2}
\volumeyear{1999}
\volumename{Proceedings of the Kirbyfest}
\pagenumbers{343}{348}
\papernumber{19}
\received{30 December 1998}\revised{12 April 1999}
\published{20 November 1999}

\newtheorem{thm}{Theorem}

\def\R{\mathbb{R}}

\title{Cusp equivalence between smooth embeddings\\of the 2--sphere in 4--space}
\shorttitle{Cusp equivalence}
\asciititle{Cusp equivalence between smooth embeddings of the 2-sphere in 4-space}

\author{Takao Matumoto}

\address{Department of Mathematics, Faculty of Science, Hiroshima University
\\Higashi-Hiroshima 739-8526, Japan }
\email{matumoto@math.sci.hiroshima-u.ac.jp}
\primaryclass{57Q45}\secondaryclass{57R45}

\keywords{Cusp, generic deformation of maps, smooth 2--knots, 
stable equivalence}
\asciikeywords{Cusp, generic deformation of maps, smooth 2-knots, 
stable equivalence}

\begin{document}

\begin{abstract}
If the fundamental group of the complement of a smooth embedding $f\co S^2 \subset \R^4$ is a cyclic group,
the map can be deformed to the standard embedding by a generic one-parameter family with
at most cusp singularities.
If two smooth embeddings are connected by such a deformation, they will be called cusp equivalent.
We will discuss the relation of three equivalences of smooth 2--knots $S^2 \subset \R^4$; cusp equivalence,
stable equivalence and weakly stable equivalence. \end{abstract}

\asciiabstract{If the fundamental group of the complement of a smooth 
embedding f: S^2 \subset R^4 is a cyclic group, the map can be
deformed to the standard embedding by a generic one-parameter family
with at most cusp singularities.  If two smooth embeddings are
connected by such a deformation, they will be called cusp equivalent.
We will discuss the relation of three equivalences of smooth 2-knots
S^2 \subset R^4; cusp equivalence, stable equivalence and weakly
stable equivalence.}

\makeshorttitle

\section{Introduction}

First we present the following theorem. We say that $f_\lambda$ has a {\it cusp}
singularity at $\lambda=0$
if $f_\lambda (x, y) = (x^2, y, x(\lambda - x^2 - y), \pm xy)$ or $f_\lambda (x, y) = (x^2, y, x(- \lambda - x^2 - y), \pm xy)$ for some local coordinates of $S^2$ and $\R^4$ which are compatible with the orientations. The former is called a {\it cusp birth} and the latter a {\it cusp death}. 

\begin{thm}
If the fundamental group of the complement of a smooth embedding $f\co S^2 \subset \R^4$ is a cyclic group, then there is a generic one-parameter family of smooth maps $f_t\co S^2 \rightarrow \R^4$ such that
$f_0$ is the standard embedding,
$f_1=f$ and $f_t$ is a self-transverse immersion except for finitely many $t$'s where $f_t$ has only one cusp singularity.
\end{thm}

If $f_t$ has a first cusp singularity at $t_1$, then $f_t$ is an embedding for $t<t_1$ and $f_t$
has a self-intersection point for $t_1<t<t_1+ \epsilon $.
The second cusp either leads a birth of another self-intersection point or leads to the death of the first self-intersection point and so on. A self-intersection point appears or disappears at a cusp birth or a cusp death respectively. Theorem~1 is easily extended to a
connected embedded surface in $\R^4$.
In the case of non-orientable surfaces there are several unknotted embeddings which are not isotopic because they have different normal Euler numbers but they are connected by such a deformation. Kamada also gave a proof of the generalized version of Theorem~1 in [2].

Let $f_0(S^2)$ and $f_1(S^2)$ be the image of two embedded 2--spheres. These 2--knots can be equivalent in the following three senses. 

(1)\qua They are {\it cusp equivalent} if they are connected by a generic one-parameter family $\{f_t\}$ of smooth maps such that $f_t$ is a self-transverse immersion except for finitely many values of $t$ at which $f_t$ has only one cusp singularity. 

(2)\qua  They are {\it stably equivalent} if they are ambient isotopic after trivial 1--handles are attached. That is, they are ambient isotopic after a genus $g$ unknotted orientable surface is attached as a connected summand to each.

(3)\qua  They are {\it weakly stably equivalent} if, for some $n \geq 0$, they are ambient isotopic in $\R^4 \# (\#_n S^2 \times S^2)$. 

If $f_t$ is a generic deformation with at most cusp singularities, then the fundamental group of the complement is kept constant during the deformation. We have the following theorem.

\begin{thm}
Assume that two smooth 2--knots $f_0\co S^2 \subset \R^4$ and $f_1\co S^2 \subset \R^4$
are cusp equivalent. Then they are stably equivalent. \end{thm}

Theorem 2 is easily generalized to the case of embedded orientable surfaces of positive genus. In the case of embedded non-orientable surfaces we need the additional condition that they have the same normal Euler numbers.

\begin{thm}
Stably equivalent 2--knots
$S^2_0 \subset \R^4$ and $S^2_1 \subset \R^4$ are weakly stably equivalent.
\end{thm}

Theorem~1 should be known to specialists who are familiar both with
2--knot theory and singularity theory.
But since it would give a starting point for a possible proof of the
smooth 4--dimensional unknotting conjecture, I would like to dedicate this short paper to Professor R Kirby.
I would also like to thank the referee who read the manuscript carefully
and offered a beautiful figure to simplify the proof of Theorem~2. 

\section{Proof of the Theorems}

\begin{proof}[Proof of Theorem 1]
Since any two maps of $S^2$ into $\R^4$ are mutually homotopic,
there is a generic one-parameter family of smooth maps $f_t\co  S^2 \to \R^4$ connecting $f_0$ and $f_1$.
By dimensional reasoning, the
generic births and deaths of self-intersection points have only two types; finger moves and cusps.
In fact at a finite number of $t$'s the image of $f_t$ has a cusp or a partial tangency of two local sheets.
By a finger move we mean the local deformation including the latter
non-transverse case; a finger pushes a local sheet until it
penetrates another local sheet.
We get a pair of self-intersection points after the penetration with one plus and one minus intersection number.
Note that the reverse process of the finger move is a Whitney trick.
We may think that the finger move follows along a curve $\alpha$ that connects two points of $S^2$.
If the fundamental group of the complement is a cyclic group generated by the element encircling the surface,
the curve $\alpha$ is homotopically trivial relative boundary and hence isotopically trivial;
so, we get a 2--disk $D$ with two corners and $\partial D= \alpha \cup (D \cap S^2)$.
Now it is not difficult to decompose this finger move into two cusp births.

The reverse deformation around the cusp birth, that is, the cusp death process can be described by
the collapsing of a 2--disk $D$ with one corner such that $\partial D = D \cap \{ {\text{\rm immersed }}S^2 \}$ and the corner point on the boundary is the self-intersection point.
Moreover we have two types of cusps for the cusp birth;
one gives a positive intersection point and another gives a
negative one for the appearing immersed surface.
If the finger move is isotopically trivial, we can construct a pair of disjoint collapsing 2--disks for the pair of positive and negative self-intersection points and get a pair of positive and negative cusps. This completes the proof of Theorem~1. For another proof one may consult the manuscript written by Kamada [2]. \end{proof}

Before giving a proof of Theorem~2
we recall the definition of trivial 1--handles. An embedded product, $B \times I$, of a 2--disk $B$ and an interval $I$ in $\R^4$ is said to {\it span} $S^2$ {\it as a 1--handle} if $(B \times I)\cap S^2 = B \times \{0,1\} \cap S^2$ and the surface obtained by the surgery
according to $B\times I$ has an
orientation compatible with the original one. The circle $\beta=\partial B \times 1/2 $ parallel to the boundary of $\partial B \times I$ is considered. The 1--handle is called {\it trivial}
if there is a 2--disk $D$ such that
$D \cap ((B \times I) \cup S^2) = \partial D$ and $\partial D \cap \beta =$ a point.

\begin{figure}[ht!]
\vglue 0.1in
\cl{\includegraphics[scale=0.4]{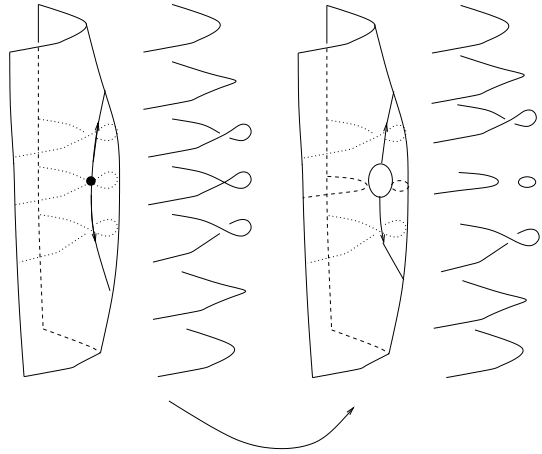}} 
\vglue -0.1in
\caption{Smoothing a double point}
\label{fig1}
\end{figure}

By a {\it smoothing of a transverse double point} we mean replacing a neighborhood
of the transverse double point with a standard annulus connecting the components of a Hopf link in the boundary 3--sphere. When the surface is oriented, the smoothing should respect the orientation and the isotopy class of the resulting surface is uniquely determined. Since a cusp birth is described also by
$f_\lambda (x, y) = (x^2, y, x(\lambda - x^2 - y^2), \pm xy)$, an
oriented smoothing of the double point appearing at the cusp birth is shown in Figure~\ref{fig1}.

We can find easily $B \times I$ and $D$ of a trivial 1--handle in the figure. If the sign of the cusp birth is reversed, the upper and lower curves at the crossing in the figure are reversed. So, we obtain:

\newtheorem*{nonumberlemma}{Lemma}
\begin{nonumberlemma} Let $F_0$ be an oriented immersed surface in $\R^4$, and let $F_1$ be a surface obtained from $F_0$ by a single cusp birth.
Let $F_2$ be obtained from $F_1$ by smoothing the double point that appears following the cusp birth.
Then $F_2$ is obtained from $F_0$ by surgery along a trivial 1--handle.
\end{nonumberlemma}

\begin{proof}[Proof of Theorem 2]
Let $f_1(S^2)=S^2_1$ and $f_2(S^2)= S^2_2$ be embedded 2--spheres in $\R^4$ that are cusp equivalent.
The one-parameter family $\{f_t\}$ can be deformed so that any cusp birth takes place at $t <1/2$
and any cusp death takes place at $ t > 1/2$. So we may assume that there is an
immersed 2--sphere $f_{1/2}(S^2)=S$ in $\R^4$ such that $S$ can be obtained from either $S^2_1$ or $S^2_2$ by attaching $m$ cusps for some $m$. Let $F$ be the embedded surface obtained from $S$ by smoothing all the double points of $S$. The result follows from the previous lemma. \end{proof}

I thank the referee for offering the use of Figure~\ref{fig1} which simplifies the proof of Theorem~2 greatly.

\begin{proof}[Proof of Theorem 3]
The surface which is obtained by attaching $n$ trivial 1--handles is a connected sum of the original 2--knot and an
unknotted surface $F$ of genus $n$ in the sense of [1].
Since the bounding solid tori for the unknotted surfaces $F$ are ambient isotopic to each other by Corollary~1.6 of [1],
we may assume that there is an orientation preserving diffeomorphism $h$ of $\R^4$ to itself
which not only satisfies $h(S^2_0 \# F)= S^2_1 \# F$ but also preserves the circles $\beta$'s determined by the trivial 1--handles.
Now we consider each trivial 1--handle.
Perform surgery along the loop on the embedded surface which intersects the circle
$\beta$ at one point for the trivial 1--handle recursively. We have to respect the stable framing when making the surgery.
Then we get a 2--knot in the new ambient manifold, the connected sum of $\R^4$ with $n$
copies of $S^2 \times S^2$. Let $D_s$ be the core 2--disk of the surgery disk and $D_o$ the original Whitney type 2--disk for the trivial 1--handle. The 2--disks $D_o$'s depend on each of the original 2--knots.
If we do surgery on the new ambient manifold along each 2--sphere $D_s\cup D_o$,
we get the original 2--knot in $\R^4$. So, doing the reversing surgery $n$ times for the complement of
each of $S^2_0$ and $S^2_1$ again, we get Theorem~3. \end{proof}

\section{On the generalization}

The generalization of Theorem~1 to
an embedded connected surface in $\R^4$ is straightforward.

In the case of $\R P^2$ we have a generic one-parameter family with only a pair of cusp
births and deaths connecting the unknotted surfaces with normal Euler numbers $2$
and $-2$.
Here the normal Euler number means half of the intersection number of the zero section of the induced disk bundle for the normal bundle of the immersed unoriented surface with respect to its orientable 2--covering.

The proof of Theorem~2 is easily generalized to an oriented surface of positive genus because the sign of a
cusp birth or death is determined by the intersection number of the self-intersection point of the central oriented immersed surface. Also, if the intersection number of the self-intersection point is positive or negative,
its contribution to the normal Euler number of the immersed surface is $-2$ or $+2$ respectively [3] and the procedure of Lemma annuls this contribution in the case of oriented smoothing. 

But in the case of
non-orientable surfaces we cannot distinguish a positive or negative cusp after it turns into a self-intersection point. So, a double point which appears at a positive cusp birth may terminate at a negative cusp death.
In this case the smoothing double point cannot respect the local orientation near the cusp and if the cusp birth is smoothed by attaching a trivial 1--handle then the cusp death is smoothed by attaching two unknotted $\R P^2$'s with the same Euler number.
Nevertheless, if the normal Euler numbers of two embeddings $f_0$ and $f_1$ are the same, the number of positive cusp births and positive cusp deaths should be equal and similarly the number of negative ones.
Hence, if a self-intersection point
born at a positive cusp dies at a negative cusp, there is also a self-intersection point born at a negative cusp and dieing at a positive cusp and vice versa. So, if the cusp births are smoothed by attaching 1--handles, then the cusp deaths are smoothed by attaching the connected sum of unknotted tori and $\R P^2$'s with zero total Euler number. Since the connected sums of the unknotted $\R P^2$ with the unknotted torus and with two unknotted $\R P^2$'s with Euler number $-2$ and $2$ are isotopic, cusp equivalent non-orientable surfaces with the same Euler number are stably equivalent
after possibly taking connected sum of each with an unknotted $\R P^2$.

If the notion of stable equivalence for non-orientable surfaces is enlarged to include the connected sum with
unknotted $\R P^2$'s, then we can drop the hypothesis that the normal Euler numbers are the same. 

As for Theorem 3, if the 1--handles are preserved by an ambient isotopy for stably equivalent embedded surfaces, they are easily shown to be weakly stably equivalent. We expect this
to be true in general but we might have some difficulties even in the orientable case.

\Addresses\recd

\end{document}